\documentclass[a4paper]{amsart}

\usepackage[english]{babel}
\usepackage{amsmath}
\usepackage{amssymb}
\usepackage{amsthm}
\usepackage{mathtools}
\usepackage{dsfont}
\usepackage{microtype}
\usepackage{rotating}
\usepackage{graphicx}
\usepackage{tikz}
\usetikzlibrary{arrows,patterns}
\usepackage[a4paper,includeheadfoot,margin=2.54cm]{geometry}
\usepackage{xspace}

\usepackage[autostyle]{csquotes}

\usepackage{pgfplots}

\usepackage[colorlinks, bookmarksdepth=1]{hyperref}

\newcommand{\rleft}{\mathopen{}\mathclose\bgroup\left}
\newcommand{\rright}{\aftergroup\egroup\right}

\newcommand{\Nd}{{\mathds{N}}}
\newcommand{\Rd}{{\mathds{R}}}
\newcommand{\Zd}{{\mathds{Z}}}

\newcommand{\N}{\Nd}
\newcommand{\R}{\Rd}
\newcommand{\Z}{\Zd}

\newcommand{\Pm}{{\mathcal{P}}}
\newcommand{\Tm}{{\mathcal{T}}}

\newtheorem{thm}{Theorem}[section]

\theoremstyle{definition}
\newtheorem{defn}[thm]{Definition}
\newtheorem{remark}[thm]{Remark}

\newcommand{\set}[1]{\rleft\{ {#1} \rright\}}
\newcommand{\with}{\colon}
\newcommand{\ie}{i.e.,}

\subjclass[2010]{Primary: 52B20; Secondary: 52C05, 14M25}
\keywords{Lattice triangles, Ehrhart polynomial, $h^\ast$-vector, toric surfaces, sectional genus, Scott's inequality}

\begin{document}
\selectlanguage{english}

\begin{abstract}
  The Ehrhart polynomial of a lattice polygon $P$ is completely determined by the
  pair $(b(P),i(P))$ where $b(P)$ equals the number of lattice points on the boundary and
  $i(P)$ equals the number of interior lattice points. All possible pairs $(b(P),i(P))$ are completely described by a theorem due to Scott. In this note, we describe the
  shape of the set of pairs $(b(T),i(T))$ for lattice triangles $T$ by finding infinitely many new Scott-type inequalities.
\end{abstract}

\title{On Ehrhart Polynomials of Lattice Triangles}
\author{Johannes Hofscheier}
\address[J.~Hofscheier and B.~Nill]{Institut f\"ur Algebra und Geometrie (IAG), Universit\"at
  Magdeburg, Geb\"aude 03, Universit\"atsplatz 2, 39106 Magdeburg,
  Germany}
\curraddr{}
\email{johannes.hofscheier@ovgu.de, benjamin.nill@ovgu.de}
\thanks{}
\author{Benjamin Nill}
\author{Dennis \"Oberg}
\address[D. \"Oberg]{Department of Mathematics, Stockholm university (SU), 10691 Stockholm, Sweden}
\email{dennis.oberg@math.su.se}
\maketitle{}

\begin{figure}[!ht]
  \centering
  \includegraphics[width=\textwidth]{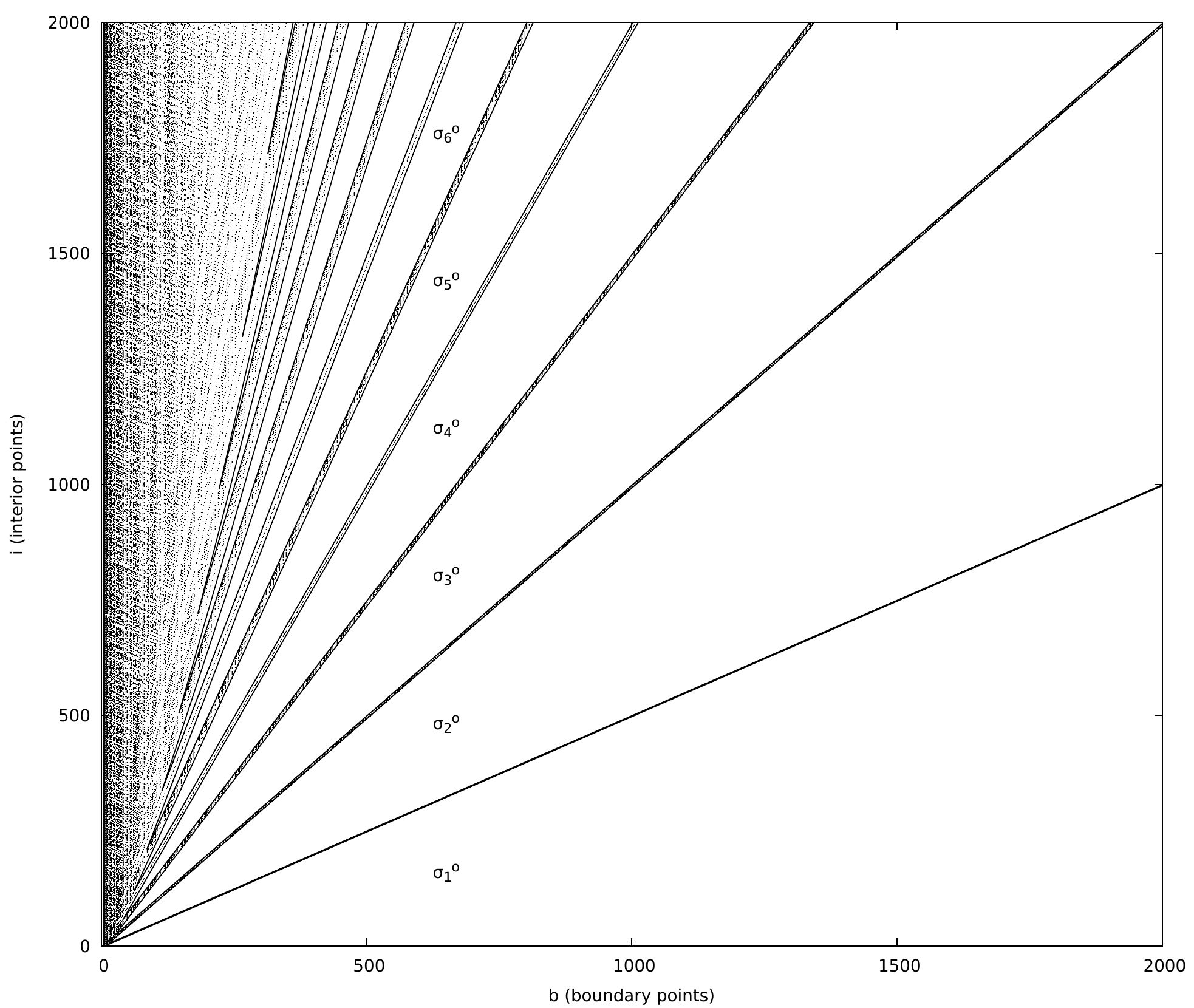}
  \caption{The points $(b(T),i(T))$ for lattice triangles $T$ together with the
    open cones $\sigma_c^\circ$.}
  \label{fig:bi2000}
\end{figure}

\section{Introduction}
\label{sec:introduction}

A \emph{lattice polygon} $P \subseteq \R^2$ is the two-dimensional convex hull
of finitely many \emph{lattice points}, \ie points in $\Z^2$. Two lattice
polygons are \emph{equivalent} if they are mapped onto each other by an
affine-linear automorphism of $\R^2$ which maps $\Z^2$ onto $\Z^2$. Let $b =
b(P)$ (resp.~$i = i(P)$) be the number of lattice points contained in the
boundary (resp.~in the interior) of $P$.

\smallskip

Pick's Theorem \cite{Pick} allows to compute the area $a(P)$ of a lattice
polygon $P$ from $b(P)$ and $i(P)$:
\begin{equation}
  \label{eq-pick}
  a(P)=i(P)+\frac{b(P)}{2}-1\text{.}
\end{equation}
The {\em Ehrhart polynomial} of $P$ is given by $|(k P) \cap \Z^2| = a(P)
k^2+\frac{b(P)}{2} k + 1$ (for $k \in \Z_{\ge 0}$). We refer to the textbook \cite{BR:Computing}. Therefore, the study of
Ehrhart polynomials of lattice polygons reduces to the study of the set $\Pm$ of
tuples $( b(P), i(P) )$ for lattice polygons $P$. In 1976 Scott showed the
following result:

\begin{thm}[Scott]
  \label{thm:scott}
  For a lattice polygon $P$ with $i(P) \ge 1$ either $(b(P),i(P)) = ( 9, 1 )$ or
  $b(P) \le 2i(P) + 6$ holds.
\end{thm}

As described in \cite{Haase}, this implies a complete description of Ehrhart
polynomials of lattice polygons:
\[
  \Pm = \{(b,0) \,:\, b \in \Z_{\ge 3}\} \cup \{(b,i) \,:\, 3 \le b \le 2i+6,\,
  1 \le i\} \cup \{(9,1)\} \text{.}
\]

\smallskip

\begin{figure}[!ht]
  \begin{sideways}
    \centering
    \includegraphics[height=\textwidth]{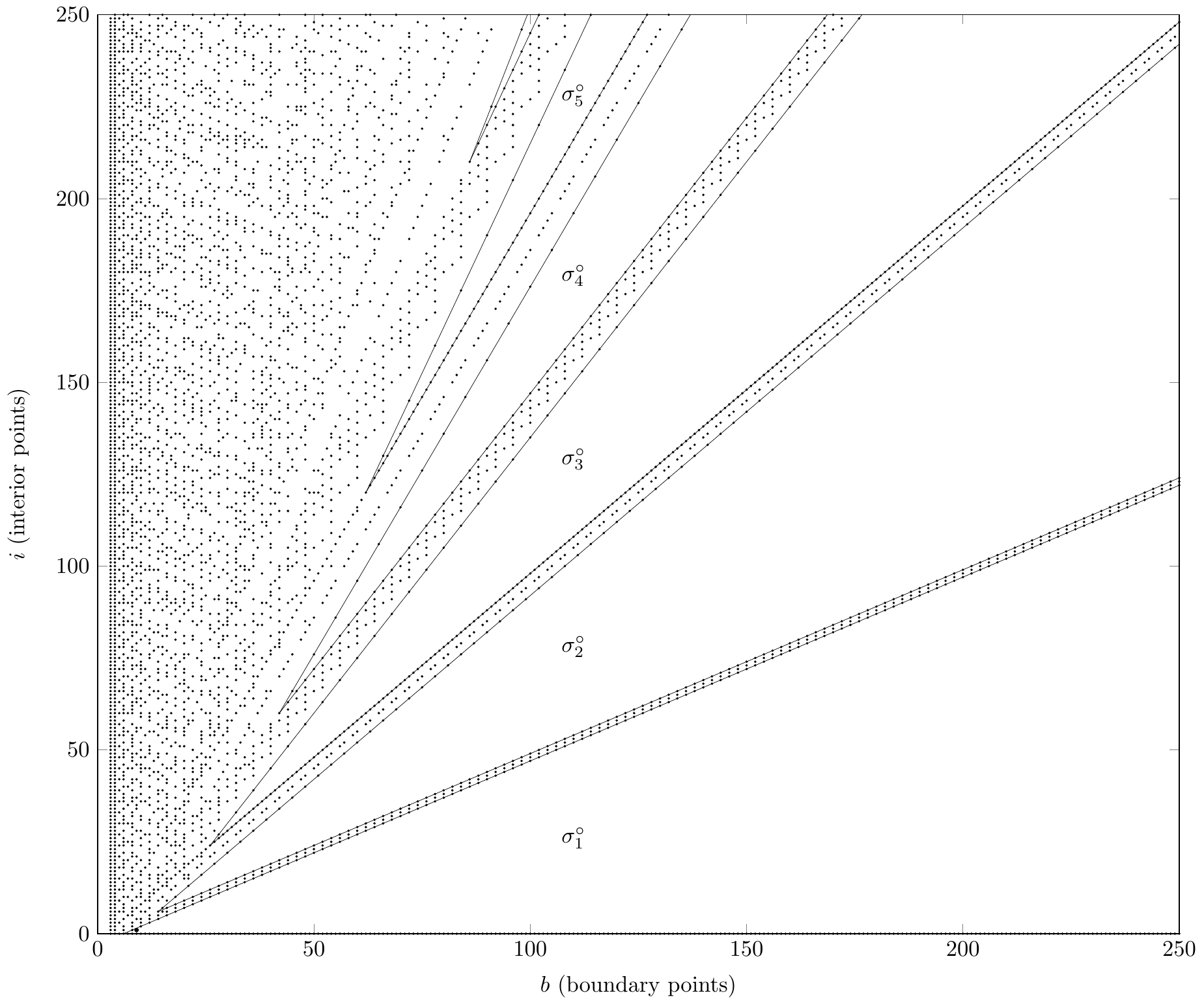}
  \end{sideways}
  \caption{Zooming into Figure \ref{fig:bi2000} for the values $0 \le b \le 250$.}
  \label{fig:bi250}
\end{figure}

In this note, we investigate the subset $\Tm \subseteq \Pm$ of tuples
$(b(T),i(T))$ for \emph{lattice triangles} $T \subseteq \R^2$. Since in Ehrhart
theory results are often reduced to the case of lattice simplices, we are
interested in understanding what this reduction means for the set of Ehrhart
polynomials in the simplest case of dimension two. Rather surprisingly, it turns
out that the structure of the set $\Tm$ is richer than one might have expected,
as the reader can see in Figure~\ref{fig:bi2000} and
Figure~\ref{fig:bi250}. While the structure of $\Tm$ is still to be fully
understood, we explain here the appearance of the conspicuous \enquote{spikes}.

For this, let us introduce the following open affine cones (see Figure
\ref{fig:bi2000} and Figure \ref{fig:bi250}).
\begin{defn}
  We set $\sigma_c^\circ \coloneqq \set{ (b,i) \in \R_{\ge0}^2 \with
    \tfrac{c-1}{2} b - ( c - 1 ) < i < \tfrac{c}{2} b - c ( c + 2 ) } \subseteq
  \R^2$ for $c \in \Z_{\ge 1}$.
\end{defn}

It is straightforward to check that the closures of these cones are pairwise
disjoint. In this notation, Scott's theorem (Theorem~\ref{thm:scott}) is
equivalent to the statement $\rleft( \Pm \setminus \{ ( 9, 1 ) \} \rright) \cap
\sigma_1^\circ = \emptyset$.  As suggested by Figure~\ref{fig:bi2000}, the
complement of $\Tm$ contains infinitely many more components.

\begin{thm}
  \label{thm:ltc}
  We have $\Tm \cap \sigma_c^\circ = \emptyset$ for $c \in \Z_{\ge 2}$.
\end{thm}

For $c \in \Z_{\ge1}$ let $\widetilde{\sigma}_c$ be the translate of the closure
of $\sigma^\circ_c$ so its apex is at the origin. As the cones
$\widetilde{\sigma}_c$ cover the positive orthant, we see that there are no
two-dimensional open affine cones in $\R^2_{\ge0}$ that are disjoint from all of
the cones $\sigma_c^\circ$.

\begin{remark}
  \leavevmode
  \begin{enumerate}
  \item There is a purely number-theoretic criterion to check whether
    a given pair $(b,i)$ is in $\Tm$. We have $(b,i) \in \Tm$ if and
    only if there exist integers $A,B,C \in \Z$ with $A >0$ and
    $0 \le B < C$ such that $b = A + \gcd( B, C ) + \gcd( B - A, C )$
    and $i = (AC - b) / 2 + 1$. In this case, the triangle with
    vertices $(0,0),(A,0),(B,C)$ can be chosen. These statements
    follow easily by considering Hermite normal forms (for details,
    see \cite{Bachelor}).
  \item Let us note that for $c \ge 1$, the apex of the closure of the
    cone $\sigma_c^\circ$ is $( 2c^2 + 2c + 2, c^3 - c ) \in \Tm$,
    realized by the lattice triangle with vertices $( 0, 0)$,$(2c^2+2c, 0 )$,$(1,
    c )$. Moreover, every pair $(b,i) \in \Z^2$ on the lower facet of the
    closure of the cone $\sigma_c^\circ$ lies in $\Tm$ and is realized
    by the lattice triangle with vertices
    $(0,0)$,$\rleft( \rleft( 2i+b-2 \rright) / c, 0 \rright)$,$(B,c)$
    for $B \in \Z_{\ge0}$ with
    $\gcd(B,c) = \gcd \rleft( B - \rleft( 2i + b - 2 \rright) / c, c
    \rright) = 1$.
    We also find infinitely many pairs
    $(k(c+1),kc(c+1)/2-c(c+2)) \in \Tm$ (for $k \in \Z_{\ge 2c+1}$) on
    the upper face of the closure of the cone $\sigma_c^\circ$,
    realized by the lattice triangles with vertices
    $(0,0)$,$((k-2)(c+1),0)$,$(0,c+1)$. These statements follow from
    elementary number-theoretic considerations.
  \item The reader might notice many missing lines of slope $-1/2$ in
    Figure~\ref{fig:bi250}. The reason is that there are only two lattice
    triangles of prime normalized volume $2a$. More precisely, for odd primes
    $p$ we have $\Tm \cap \set{ i = -\tfrac{1}{2}b + \tfrac{p}{2} + 1 } =
    \set{(3,\tfrac{p-1}{2}), (p+2,0)}$, for details see \cite{Bachelor}. This
    follows also from Higashitani's study of lattice simplices with prime
    normalized volume \cite{Hig}.
  \item From the pictures, it is visible that the points of $\Tm$ in the
    \enquote{spikes} form periodic patterns. It seems to be an interesting open question to
    make this observation precise.
  \item Scott's inequality (Theorem~\ref{thm:scott}) follows from inequalities by del
    Pezzo and Jung \cite{Jung,delPezzo} (see also \cite{Schicho}): any
    rational projective surface with degree $d$ and sectional
    genus $p > 0$ satisfies $d \le 4p+4$ if $(d,p) \not= (9,1)$. Now,
    Theorem~\ref{thm:ltc} can be translated into algebraic geometry as follows:
    there exists no toric projective surface with Picard number one, degree $d$,
    and sectional genus $p$ that satisfies
    \[
      \tfrac{2(c+1)}{c} (c+p) < d < \tfrac{2c}{c-1} p
    \]
    for an integer $c \ge 2$. It would be interesting to see whether this is a
    special case of a more general algebro-geometric statement.
  \end{enumerate}
\end{remark}

We remark that in an upcoming paper of the first author and Higashitani Theorem
\ref{thm:ltc} will be used as the base case of a generalization for lattice
simplices of dimension greater than two.

\section*{Acknowledgments}
This project was originally inspired by an email correspondence with Tyrrell
McAllister. The authors would also like to thank Gabriele Balletti for his
supply of computational data, e.g., Figure~\ref{fig:bi2000} and
Figure~\ref{fig:bi250} is due to him, as well as many fruitful
discussions. Theorem~\ref{thm:ltc} proves a conjecture in the bachelor's thesis
of the third author \cite{Bachelor}, advised by the second author at Stockholm
University.  The second author is an affiliated researcher with Stockholm
University and partially supported by the Vetenskapsr{\aa}det
grant~NT:2014-3991.

\section{Proof of Theorem \ref{thm:ltc}}
\label{sec:proofs}

Our proof uses the ideas of Scott's original proof of
Theorem~\ref{thm:scott}. For the convenience of the reader we will give complete
arguments without assuming prior knowledge of \cite{Scott}. Let $T$ be a lattice
triangle with area $a$, number of boundary lattice points $b$, and number of
interior lattice points $i$. We assume that $T$ satisfies for some $c \ge 2$
the inequalities
\begin{equation}
  \label{eq:zz}
  c \rleft( \tfrac{b}{2} - 1 \rright) < a < \tfrac{b}{2} (c+1) -
  (c+1)^2 \text{.}
\end{equation}
We will show that this situation cannot exist. 

\smallskip

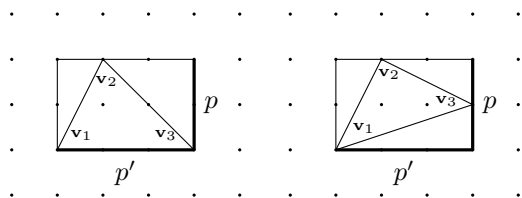
\begin{figure}[!ht]
  \centering
  \begin{tikzpicture}[scale=0.6]
    \clip (-1.04, -1.04) -- (4.04, -1.04) -- (4.04, 3.04) -- (-1.04, 3.04) -- cycle;
    \draw (0,0) -- (3,0)  -- (3,2) -- (0,2) -- cycle;
    \draw (0,0) node[above right,xshift=0.5mm] {\tiny $\mathbf{v}_{1}$} -- (1,2)
    node[below=1mm,xshift=0.5mm] {\tiny $\mathbf{v}_{2}$} -- (3,0)
    node[above left,xshift=-1mm] {\tiny $\mathbf{v}_{3}$} -- cycle;
    \draw[very thick] (3,0) -- (3,2) node[midway,right]{$p$};
    \draw[very thick] (0,0) -- (3,0) node[midway,below]{$p'$};
    \foreach \x in {-1,...,4} \foreach \y in {-1,...,3} \fill (\x, \y) circle (1pt);
  \end{tikzpicture}\hspace*{0.5cm}
  \begin{tikzpicture}[scale=0.6]
    \clip (-1.04, -1.04) -- (4.04, -1.04) -- (4.04, 3.04) -- (-1.04, 3.04) -- cycle;
    \draw (0,0) -- (3,0) -- (3,2) -- (0,2) -- cycle;
    \draw (0,0) node[above right=1mm] {\tiny $\mathbf{v}_{1}$} -- (3,1)
    node[left=0.75mm,yshift=0.5mm] {\tiny $\mathbf{v}_{3}$} -- (1,2)
    node[below,xshift=1mm] {\tiny $\mathbf{v}_{2}$} -- cycle;
    \draw[very thick] (3,0) -- (3,2) node[midway,right]{$p$};
    \draw[very thick] (0,0) -- (3,0) node[midway,below]{$p'$};
    \foreach \x in {-1,...,4} \foreach \y in {-1,...,3} \fill (\x, \y) circle (1pt);
  \end{tikzpicture}
  \caption{Two examples of bounding boxes used in the proof of Theorem \ref{thm:ltc}.}
  \label{fig:ltc-bb}
\end{figure}

By replacing $T$ with an equivalent lattice triangle, we may assume that $T$ is contained in a bounding box (\ie a rectangle whose edges are parallel to the coordinate axes and that is minimal with respect to the inclusion of $T$) with vertical side length $p$ such that $p$ is minimal among all such choices. For an illustration, see Figure \ref{fig:ltc-bb}. Let us note that $p$ equals the lattice width of $T$. We denote the horizontal side length of the bounding box by $p'$. We observe that necessarily $p' \ge p$ since switching coordinates yields an equivalent triangle.

\smallskip
  
Let us prove 
\begin{equation}
  \label{eq:a-lb-p}
  a \ge \tfrac{p^2}{4}
\end{equation}

For this, we denote by $\delta$ the minimal distance of the $x$-coordinate of a vertex of $T$ (denoted by $\mathbf{v}_2$) on the top edge of the rectangle from the $x$-coordinate of a vertex of $T$ (denoted by $\mathbf{v}_1$) on the bottom edge of the rectangle. By an integral, unimodular shear leaving the horizontal line through the bottom edge of the rectangle invariant, we can achieve $\delta \le
  \tfrac{p}{2}$. By possibly flipping along the horizontal or vertical axis, we may also assume that $\mathbf{v}_2$ has $x$-coordinate greater than or equal to that of $\mathbf{v}_1$, and the third vertex of $T$ (denoted by $\mathbf{v}_3$) has $x$-coordinate greater than or equal to that of $\mathbf{v}_2$ (recall that $p' \ge p > \delta$).
  
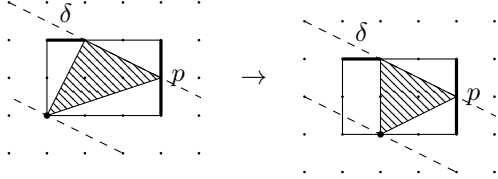
\begin{figure}[!ht]
  \centering
  \parbox{2.6cm}{
    \begin{tikzpicture}[scale=0.5]
      \clip (-1.04, -1.04) -- (4.04, -1.04) -- (4.04, 3.04) -- (-1.04, 3.04) -- cycle;
      \draw (0,0) -- (3,0) -- (3,2) -- (0,2) -- cycle;
      \draw[pattern=north west lines] (0,0) -- (3,1) -- (1,2) -- cycle;
      \draw[very thick] (3,0) -- (3,2) node[midway,right]{$p$};
      \draw[very thick] (0,2) -- (1,2) node[midway,above=1mm]{$\delta$};
      \filldraw (0,0) circle (2pt);
      \draw[dashed] (-1,3) -- (5,0);
      \draw[dashed] (-2,1) -- (2,-1);
      \foreach \x in {-1,...,4} \foreach \y in {-1,...,3} \fill (\x, \y) circle (1pt);
    \end{tikzpicture}
  } $\quad \rightarrow \quad$
  \parbox{2.6cm}{
    \begin{tikzpicture}[scale=0.5]
      \clip (-1.04, -1.04) -- (4.04, -1.04) -- (4.04, 3.04) -- (-1.04, 4.04) -- cycle;
      \draw (0,0) -- (3,0) -- (3,2) -- (0,2) -- cycle;
      \draw[pattern=north west lines] (1,0) -- (3,1) -- (1,2) -- cycle;
      \draw[dashed] (-1,3) -- (5,0);
      \draw[dashed] (-1,1) -- (3,-1);
      \filldraw (1,0) circle (2pt);
      \draw[very thick] (3,0) -- (3,2) node[midway,right]{$p$};
      \draw[very thick] (0,2) -- (1,2) node[midway,above=1mm]{$\delta$};
      \foreach \x in {-1,...,4} \foreach \y in {-1,...,3} \fill (\x, \y) circle (1pt);
    \end{tikzpicture}
  }
  \caption{Illustration of the proof of inequality (\ref{eq:a-lb-p}).}
  \label{fig:a-lb-p}
\end{figure}

  Now, we move the bottom vertex $\mathbf{v}_1$ horizontally to the right until it has the same $x$-coordinate as that of $\mathbf{v}_2$. We observe that the area of the obtained triangle is bounded by
  the area of $T$ (see Figure \ref{fig:a-lb-p}), \ie
  \[
    a \ge \tfrac{p(p'-\delta)}{2} \ge \tfrac{p^2}{4}
 \]
 where for the second inequality we used $\delta \le \tfrac{p}{2}$ and $p' \ge p$. This finishes the proof of \eqref{eq:a-lb-p}.
 
\smallskip
          
We may assume that there is only one vertex of $T$ on the top edge of the bounding box (otherwise, flip horizontally). Let us denote by $q$ the length of the intersection of $T$ with the bottom edge of the bounding box, so $q=0$ if and only if there is only one vertex of $T$ on the bottom edge.

As each horizontal line between $y=0$ and $y=p$ cuts the boundary of $T$ in two
points (see Figure \ref{fig:ltc-bb}), we obtain
\begin{equation}
  \label{eq:ab-lb-pq}
  b \le q + 2p \text{.}
\end{equation}
We reformulate the inequality on the right hand side of (\ref{eq:zz}) as
\[
  \tfrac{2a}{c+1} + 2(c+1) < b\text{.}
\]
By Pick's Theorem \eqref{eq-pick}, we have $a \in \tfrac{1}{2} \N$, and thus the strict
inequality becomes $\tfrac{2a+1+2(c+1)^2}{c+1} \le b$. Combining this with
(\ref{eq:ab-lb-pq}) gives
\begin{equation}
  \label{eq:main}
  \tfrac{2a+1+2(c+1)^2}{c+1} \le b \le q + 2p\text{.}
\end{equation}
  
Let us assume $q=0$. Plugging in (\ref{eq:a-lb-p}) yields
\[
  \tfrac{p^2/2+1+2(c+1)^2}{c+1} \le 2p \Rightarrow p^2 - 4(c+1)p + 2 + 4(c+1)^2 \le 0\text{.}
\]
This is a contradiction as the discriminant of this quadratic polynomial in $p$
is negative.

\smallskip

Hence, $q > 0$, so  $a = \tfrac{pq}{2}$. Plugging this into
\eqref{eq:main} implies
\[
\tfrac{pq+1+2(c+1)^2}{c+1} \le q + 2p\text{.}
\]
Solving for $q$ yields
\begin{equation}
  \label{eq:q}
  (p-c-1)q \le 2(c+1)(p-c-1)-1\text{.}
\end{equation}
We see that $p \not= c+1$.
  
Let us assume $p < c + 1$, \ie $p \le c$. Clearly, $b \ge q+2$. We plug this in
the inequality on the left hand side of (\ref{eq:zz}) and get a contradiction, namely
\[
  \tfrac{cq}{2} < a =
  \tfrac{pq}{2} \le \tfrac{cq}{2} \text{.}
\]

Hence, we have $p > c + 1$. We deduce from (\ref{eq:q}) 
\begin{equation}
  \label{new}
  q < 2 (c+1) \text{.}
\end{equation}

Let us translate the left bottom vertex of $T$ into the origin. By applying an integral, unimodular shear leaving the horizontal line through the bottom edge of the rectangle invariant, we can get an equivalent triangle such that the $x$-coordinate of the top vertex of $T$ is in $[0, p)$. As $p' \ge p$, this implies $q=p'$, e.g., as in the left example of Figure \ref{fig:ltc-bb}.

The inequality on the left hand side of (\ref{eq:zz}) is equivalent to $b <
\tfrac{2a}{c}+2$. As $a = pq/2 \in \tfrac{1}{2} \N$, the strict inequality becomes
\[   
b \le \tfrac{2a - 1}{c} + 2 = \tfrac{pq-1}{c} + 2 \text{.}
\]
Combining this with the inequality on the left hand side of (\ref{eq:main}), we
obtain
\[
  \tfrac{pq+1+2(c+1)^2}{c+1} \le b \le \tfrac{pq-1}{c} + 2\text{.}
\]
Solving for $pq$ yields $2c(c^2+c+1) + 1 \le pq$. As $q=p' \ge p$ and $p > c+1$,
the previous inequality combined with \eqref{new} implies
\begin{equation}
  \label{eq:q-sq}
  2c(c^2+c+1) < pq \le q^2 < 4(c+1)^2
  \Rightarrow 2c^3 -2c^2 -6c -4 < 0\text{.}
\end{equation}
A straightforward computation shows that this is only possible for $c \le2$. As
$c \ge 2$, we deduce $c=2$. Plugging this again into (\ref{eq:q-sq}) we obtain
$28 < q^2 < 36$, a contradiction. \hfill\qed

\bibliographystyle{amsalpha}
\bibliography{epolt}

\end{document}